\documentclass[amssymb,amsfonts,reqno]{amsart}
\usepackage{latexsym}

\usepackage{graphicx}

\newcommand{\gtrsim}{\raisebox{-0.13cm}{~\shortstack{$>$ \\[-0.07cm]
      $\sim$}}~}
      
\newcommand{\lesssim}{\raisebox{-0.13cm}{~\shortstack{$<$ \\[-0.07cm]
      $\sim$}}~}

\def\be#1{\begin{equation}\label{#1}}
\def\bas{\begin{align*}}
\def\eas{\end{align*}}
\def\bi{\begin{itemize}}
\def\ei{\end{itemize}}

\usepackage{pdfpages}
\usepackage[framemethod=TikZ]{mdframed}

\mdfdefinestyle{JimFrame}{%
    linecolor=red,
    outerlinewidth=2pt,
    roundcorner=20pt,
    innertopmargin=\baselineskip,
    innerbottommargin=\baselineskip,
    innerrightmargin=25pt,
    innerleftmargin=25pt,
    backgroundcolor=white!50!white}

\theoremstyle{plain}
   \newtheorem{theorem}[subsection]{Theorem}

   \newtheorem{lemma}[subsection]{Lemma}
   \newtheorem{corollary}[subsection]{Corollary}

\begin{document}

\author{Paolo Ciatti}
\address{Universit\`a degli Studi di Padova, 
Via Marzolo 9, 35131 Padova,
Italia}
\email{paolo.ciatti@unipd.it}

\author{James Wright}
\address{Maxwell Institute of Mathematical Sciences and the School of Mathematics, University of
Edinburgh, JCMB, The King's Buildings, Peter Guthrie Tait Road, Edinburgh, EH9 3FD, Scotland}
\email{J.R.Wright@ed.ac.uk}




\title[Strongly singular integrals on stratified groups]{
Strongly singular integrals on stratified groups}

\maketitle

\centerline{\it In honour of Fulvio Ricci on his 70th birthday}

\begin{abstract} We consider a class of spectral multipliers on stratified Lie groups 
which generalise the class of H\"ormander multipliers and include multipliers
with an oscillatory factor. Oscillating multipliers have been examined extensively in the euclidean setting
where sharp, endpoint $L^p$ estimates are well known. In the Lie group setting, corresponding
$L^p$ bounds for oscillating spectral multipliers have been established by several authors but only in 
the open range of exponents. In this paper we establish the endpoint $L^p(G)$ bound when 
$G$ is a stratified Lie group. More importantly we begin to address whether these estimates are sharp.
\end{abstract}

\section{Introduction}\label{introduction}


The following class of strongly singular convolution operators on ${\mathbb R}^n$
given by
$$
T_{a,b} f(x) \  = \ \int_{|y|\le 1} f(x- y) \, \frac{e^{i |y|^{-a}}}{|y|^b} \, dy
$$ 
where $a>0$ and $b \le n(2+a)/2$
has a rich and interesting history. In the periodic setting, they were investigated by 
Hardy who used them to construct a variety of counterexamples. 
Regarding $L^p$ boundedness properties,
Hirschman \cite{Hirsch} considered the one dimensional case and for general $n\ge 1$, 
Wainger \cite{W} established 
the sharp $L^p$ range but left open the endpoint case which C. Fefferman and Stein \cite{FS1}
accomplished using
interpolation by proving that $T_{a,n}$ is bounded on the Hardy space $H^1({\mathbb R}^n)$.
Earlier C. Fefferman \cite{F} established that $T_{a,n}$
satisfies a weak-type $(1,1)$ bound.  Chanillo \cite{Chan} extended these results
to weighted $L^p$ estimates. It is well known that when $b>n(2+a)/2$, there are
no $L^p$ estimates.

As a convolution operator, we can view $T = T_{a,b}$ as a multiplier operator 
${\widehat{Tf}}(\xi) = m(\xi) {\widehat{f}}(\xi)$ where $m = m_{\theta, \beta}$ is essentially given by
\begin{equation}\label{ss-multipliers}
m_{\theta, \beta}(\xi)  \ = \ \frac{e^{i |\xi|^{\theta}}}{|\xi|^{\theta \beta/2}}
\end{equation}
for $|\xi|$ large. Here $0<\theta = a/(1+a) < 1$ and $\beta = ((2 + a)n - 2b)/a$. We note that $m$ is bounded 
precisely when $b\le n(2+a)/2$.

The case $b = n$, or equivalently $\beta = n$ in \eqref{ss-multipliers}, corresponds to the singular integral operators
$T_{a,n}$, treated by Fefferman and Stein, whose convolution kernels just fail to be integrable. Their
multipliers $m_{\theta, n}$
are not H\"ormander multipliers but furnish
examples of multipliers with $S_{\rho,\delta}^{-m}$ symbols where $m\ge 0$ and $\rho < 1$. In
this context these multipliers were studied by H\"ormander \cite{H}. 

Note that the multipliers $m_{\theta, \beta}$ in \eqref{ss-multipliers}
with $\beta > n$ (so that $b < n$) correspond
to operators $T_{a,b}$ with integrable convolution kernels and hence are bounded
on $L^1$. 
For any $\delta > 0$, consider the analytic family $T_z^{\delta}, {\rm Re}(z) \in [0,1]$, of operators  with multipliers
\begin{equation}\label{analytic-family}
{m}_z^{\delta} (\xi) \ = \frac{e^{i|\xi|^{\theta}}}{|\xi|^{[\theta (n +\delta)/2] \, z}} \chi(\xi)
\end{equation}
where $\chi(\xi) = 0$ when $|\xi| \le 1$. Thus
$T_z^{\delta}$ is bounded on $L^2$ when $z = iy$ with $\|T_{iy}^{\delta}\|_{2\to 2}$ uniformly bounded in 
$y\in {\mathbb R}$. Also $T_z^{\delta}$ is bounded on $L^1$ when $z = 1 + iy$, again with $\|T_{1+iy}^{\delta}\|_{1\to 1}$
uniformly bounded in $y\in {\mathbb R}$.
By analytic interpolation, we see that $m_{\theta, \beta}$ is an $L^p$ multiplier in the open range
$|1/p - 1/2| <  \beta/2n$. To establish endpoint bounds, one needs to say something about
the {\it endpoint} multipliers $m_{\theta,n}$ (the case $z=1$ and $\delta = 0$ in \eqref{analytic-family}). More precisely in \cite{FS1}, Fefferman and Stein 
show that multipliers $m_{1+it}^0$ in \eqref{analytic-family} are $H^1$ multipliers with an 
operator norm at most $(1+|y|)^{n+1}$.

Fefferman and Stein developed a more general theory of multipliers which 
include the examples \eqref{ss-multipliers} as special cases.  Let $K$ be a distribution of compact support,
which is integrable away from the origin. Its Fourier transform ${\widehat K}$ is of course a function.
We make the following assumptions:
\begin{equation}\label{K-assumptions}
\begin{cases}
\int_{|x|> 2 |y|^{1-\theta}} |K(x-y) - K(x)| \, dx \ \le \ B, \ \ 0 < |y| \le 1, \\
|{\widehat K}(\xi)| \ \le \ B \, (1 + |\xi|)^{-\theta n/2} . \\
\end{cases}
\end{equation}
In \cite{FS1}, Fefferman and Stein show if $K$ satisfies \eqref{K-assumptions}, then
$|\xi|^{(n-\beta)\theta/2} {\widehat K}(\xi),  0 \le \beta < n$,
is an $L^p({\mathbb R}^n)$ multiplier when $|1/p - 1/2| \le \beta/2n$. See \cite{S1} where this
result is established in the open range $|1/p - 1/2| < \beta/2n$.  

In the papers \cite{CKW1} and \cite{CKW2}
(see also \cite{M}), Chanillo,
Kurtz and Sampson considered the cases $\theta>1$ and $\theta<0$ (here the $|\xi|$ large restriction
becomes $|\xi|$ small). Hence multipliers on ${\mathbb R}^n$ of the form
\begin{equation}\label{osc-multipliers}
m_{\theta, \beta}(\xi)  \ = \ \frac{e^{i |\xi|^{\theta}}}{|\xi|^{\theta \beta/2}} \, \chi_{\pm}(\xi)
\end{equation}
for any $\theta \in {\mathbb R}$ and $\beta \ge 0$ have been studied. Here $\chi_{+}(\xi) \equiv 0$
for $|\xi|\le 1$ when $\theta > 0$ and $\chi_{-}(\xi) \equiv 0$ when $|\xi|\ge 1$ when $\theta < 0$.

The case $\theta = 1$ is special and is related to the wave operator. The sharp range of $L^p$
bounds in this case is different from the case $\theta \not= 1$; see \cite{P} and \cite{M1}. 
We will not consider the case $\theta = 1$ and assume always $\theta \not= 1$.


In this paper we will put all these oscillating multipliers into a single, general framework 
(much like what Fefferman and Stein do in \eqref{K-assumptions} when $0 < \theta < 1$) which
strictly generalises the class of H\"ormander multipliers and furthermore we will give a unified,
purely spectral treatment which readily extends to estimates for corresponding spectral
multipliers on any stratified Lie group. 

\subsection{Notation}
 Keeping track of constants
and how they depend on the various parameters will be important for us. For the most part, constants $C$
appearing in inequalities $P \le C Q$ between positive quantities $P$ and $Q$ will be {\it absolute} or
{\it uniform} in that they can be taken to be independent of the parameters of the underlying problem. 
We will use $P \lesssim Q$ to denote $P \le C Q$ and $P \sim Q$ to denote $C^{-1} Q \le P \le C Q$.
Furthermore, we use $P \ll Q$ to denote $P \le \delta Q$ for a sufficiently small constant $\delta>0$
whose smallness will depend on the context.


{\bf Acknowledgement}: We woud like to thank Alessio Martini and Steve Wainger 
for discussing the history of the problem as well as
guiding us through the literature.

\section{The euclidean setting ${\mathbb R}^n$}\label{euclidean}
We start in the euclidean setting ${\mathbb R}^n$. 
Let $\phi \in C^{\infty}_0({\mathbb R}^n)$ be supported away from the origin and let
$m^j(\xi) := m(2^j \xi) \phi(\xi)$. It is natural to impose conditions on the $j$th pieces $m^j$.
The classical H\"ormander condition requires uniform (in $j$) control of some $L^2$ Sobolev
norm $\|m^j\|_{L^2_s}$ with $s$ derivatives. Here we want to consider not
only classical H\"ormander multipliers but also oscillating multipliers $m_{\theta, \beta}$
described in \eqref{osc-multipliers}.
Special among these are the {\it endpoint} multipliers $m_{\theta, n}$ whose bounds we interpolate with
trivial $L^2$ bounds to deduce sharp $L^p$ bounds for $m_{\theta, \beta}$ for general $\beta \ge 0$.
Hence our conditions will not only involve a smoothness parameter $s > 0$ but also an oscillation parameter
$\theta \in {\mathbb R}$ and a decay parameter $\beta\ge 0$.

For any $\theta \in {\mathbb R}$, the condition $j \theta > 0$ identifies
the frequency range of interest. In fact if $\theta > 0$, then $j \theta > 0$ corresponds to $j > 0$
or $|\xi| \ge 1$ which is the relevant frequency range indicated in \eqref{osc-multipliers}. However
if $\theta < 0$, then $j \theta > 0$ corresponds to $j < 0$ or $|\xi| \le 1$ which is the range of
interest for the oscillating multipliers in \eqref{osc-multipliers} with $\theta < 0$. Finally when
$\theta = 0$, the condition $j \theta > 0$ is vacuous.

\subsection{Our multiplier conditions} We consider the following conditions on a multiplier $m$ which will depend on parameters $s, \theta$
and $\beta$. When $j \theta \le 0$, 
we impose the standard uniform $L^2$ Sobolev norm control on the $m_j$; 
\begin{equation}\label{hypothesis-euclidean-neg}
\sup_{j:  j \theta \le 0} \| m^j \|_{L^2_s({\mathbb R}^n)} \ < \ \infty.
\end{equation}
For $j \theta > 0$, we consider the condition
\begin{equation}\label{hypothesis-euclidean-pos}
\sup_{j: j \theta >  0} \,  2^{j \theta \beta/2} \|m^j\|_{L^{\infty}({\mathbb R}^n)},  
\  \ 2^{- j \theta (2s - \beta) /2} \|m^j\|_{L^2_s({\mathbb R}^n)} \ < \ \infty. 
\end{equation}
When $\theta = 0$, the condition \eqref{hypothesis-euclidean-pos} is vacuous and \eqref{hypothesis-euclidean-neg} 
reduces to the condition
$\sup_j \|m^j\|_{L^2_s} < \infty$ and if this holds for some $s > n/2$, the classical H\"ormander
theorem states that the multiplier operator is
of weak-type $(1,1)$ and maps $H^1({\mathbb R}^n)$ boundedly into $L^1({\mathbb R}^n)$.
See \cite{S}.

One can easily verify that  the conditions \eqref{hypothesis-euclidean-neg} and  \eqref{hypothesis-euclidean-pos}
are satisfied 
for $m_{\theta, \beta}$ in \eqref{osc-multipliers} and for all $s>0$. 
Note that in \eqref{hypothesis-euclidean-pos}, the quantity $j \theta$ is always positive
and so \eqref{hypothesis-euclidean-pos} expresses a growth in the Sobolev norm $L^2_s$ of $m^j$
(when $s > \beta/2$) and a decay in the $L^2$ norm of $m^j$. 
If the condition \eqref{hypothesis-euclidean-pos} is satisfied for some $s>0$, it does not necessarily hold
for all $s' \le s$. 
Therefore we introduce $M_{\theta, \beta, s}$ consisting of those
functions $m$ which satisfiy \eqref{hypothesis-euclidean-neg} with exponent $s$ and satisfies
\eqref{hypothesis-euclidean-pos} for all exponents $s' \le s$.

\subsection{Our multiplier classes} Hence $\cup_{s>n/2} M_{0,*,s}$  
is the classical class of  H\"ormander multipliers and so  
$$
{\mathcal M}_n \ := \ \bigcup_{\theta \in {\mathbb R}\setminus\{1\}, s>n/2} M_{\theta, n, s}
$$ 
gives us a natural extension of H\"ormander multipliers. It is easy to verify that
the conditions \eqref{hypothesis-euclidean-neg} and  \eqref{hypothesis-euclidean-pos}
are independent on the choice of bump function $\phi$ and hence for any $\beta\ge 0$,
\begin{equation}\label{M-beta}
{\mathcal M}_{\beta} \ := \ \bigcup_{\theta \in {\mathbb R}\setminus\{1\}, s>n/2} M_{\theta, \beta, s} \ = \Bigl\{
|\xi|^{(n-\beta)\theta/2} m(\xi) :  m \in {\mathcal M}_n \Bigr\}.
\end{equation}
This puts us in the position to employ the analytic interpolation argument in \cite{FS1}
to deduce that $m \in {\mathcal M}_{\beta}$ is an $L^p$ multiplier in the sharp range $|1/p - 1/2| \le \beta/2n$
from $H^1$ bounds for multiplier operators associated to $m\in {\mathcal M}_n$. 

In fact one advantage of
working with ${\mathcal M}_n$ (over say, the class of multipliers arising from kernels satisfying \eqref{K-assumptions}
in the case $0<\theta < 1$) is the class ${\mathcal M}_n$ has the desirable property
that it is invariant under multiplication by $|\xi|^{iy}$ for any real $y \in {\mathbb R}$; that is, if $m \in {\mathcal M}_n$,
then $|\xi|^{iy} m(\xi)$ lies in ${\mathcal M}_n$, satisfying the bounds  \eqref{hypothesis-euclidean-neg} and \eqref{hypothesis-euclidean-pos} with polynomial growth in $|y|$. Hence for the analytic interpolation argument,
we only need to establish
that multipliers in ${\mathcal M}_n$ map $H^1$ to $L^1$ instead of showing they map $H^1$ to $H^1$ as needed in \cite{FS1}. This will be particularly useful when we move to the setting of Lie groups.

\subsection{The basic decomposition} When we analyse a multiplier $m \in {\mathcal M}_{\beta}$, we will decompose
 $m = \sum_j m_j$ where $m_j(\xi) = m(\xi) \phi(2^{-j} \xi)$ for some $\phi \in C^{\infty}_0({\mathbb R}^n)$
 supported away from the origin such that $\sum_j \phi(2^{-j} \xi) =1$ for all $\xi \not=0$. 
Note that $m^j(\xi) = m_j(2^j \xi)$ is the $j$th piece on which we impose the conditions
\eqref{hypothesis-euclidean-neg} and \eqref{hypothesis-euclidean-pos}. 
We split the multiplier $m =  m_{small} + m_{large}$ into two parts where
\begin{equation}\label{m-split}
m_{small}(\xi) \ := \ \sum_{j: j \theta \le 0} m_j (\xi) \ \ \ {\rm and} \ \ \  
m_{large}(\xi) \ := \ \sum_{j: j\theta > 0}  m_j (\xi).
\end{equation}
If $\theta = 0$, then $m = m_{small}$ and in general we note that
$m_{small}$ is a H\"ormander multipler (since \eqref{hypothesis-euclidean-neg} holds for some $s>n/2$) 
and so it is an $L^p$ multiplier for all $1<p<\infty$
(as well as a weak-type $(1,1)$ and an $H^1$ multiplier). 
We introduce the notation ${\mathcal K}_F$ to denote 
the convolution kernel associated to a multiplier $F$. 
Hence it suffices to treat the operator
$$
T^l f(x) \ = \ \sum_{j : j \theta > 0} {\mathcal K}_{m_j} * f(x) \ =: \ {\mathcal K}^{l} * f(x) 
$$ 
corresponding to the interesting frequency range where the $j$th pieces $m^j$ satisfy
\eqref{hypothesis-euclidean-pos}. 

\subsection{${\mathcal M}_n$ versus \eqref{K-assumptions} }
When $m \in M_{\theta, n, s} \subset {\mathcal M}_n$
for $0 < \theta < 1$, we claim that ${\mathcal K}^l$ 
satisfies the condition \eqref{K-assumptions} of Fefferman and Stein in \cite{FS1} (see also \cite{S1}). 
Hence for $0<\theta <1$, the class of convolution operators satisfying \eqref{K-assumptions} is larger
than the class ${\mathcal M}_n$.
In fact 
the $L^{\infty}$ condition on the $m^j$ in \eqref{hypothesis-euclidean-pos} is equivalent to the bound
$|{\widehat{{\mathcal K}^l}}(\xi)| \le B (1 + |\xi|)^{-\beta n/2}$. Furthermore we bound
$$
\int_{|x|\ge 2|y|^{1-\theta}} |{\mathcal K}^l(x - y) - {\mathcal K}^l(x) | \, dx \ \le \ \sum_{j>0} \
\int_{|x|\ge 2 |y|^{1-\theta}} |{\mathcal K}_{m_j}(x - y) - {\mathcal K}_{m_j}(x) | \, dx
$$ 
and split the sum on the right $\sum_{j\in J_1} + \sum_{j\in J_2}$ where $J_1 = \{ j>0 : 2^j \ge |y|^{-1}\}$
and $J_2 = {\mathbb N} \setminus J_1$. For the sum over $J_1$, we bound each
$$
\int_{|x|\ge 2 |y|^{1-\theta}} |{\mathcal K}_{m_j}(x - y) - {\mathcal K}_{m_j}(x) | \, dx \ \le \ 
2 \, \int_{|x| \ge |y|^{1-\theta}} |{\mathcal K}_{m_j}(x)| \, dx
$$ 
and note that if $s>n/2$,
$$
 \int_{|x| \ge |y|^{1-\theta}} |{\mathcal K}_{m_j}(x)| \, dx =
  \int_{|x| \ge 2^j |y|^{1-\theta}} |{\mathcal K}_{m^j}(x)| \, dx =
   \int_{|x| \ge 2^j |y|^{1-\theta}} |{\mathcal K}_{m^j}(x)| \, |x|^s |x|^{-s}  dx
$$
\begin{equation}\label{same-way}
\lesssim \ (2^j |y|^{1-\theta})^{-(s - n/2)} \|m^j \|_{L^2_s} \ \lesssim \ (2^{j}|y|)^{-(1-\theta)(s - n/2)}
\end{equation}
by Cauchy-Schwarz and \eqref{hypothesis-euclidean-pos}. This is summable for $j \in J_1$
leaving us to treat the sum over $J_2$. In this case we bound
\begin{equation}\label{difference-y}
\int_{|x|\ge 2 |y|^{1-\theta}} |{\mathcal K}_{m_j}(x - y) - {\mathcal K}_{m_j}(x) | \, dx \ \le \ 
|y| \, \int_{|x| \ge |y|^{1-\theta}} |\nabla{\mathcal K}_{m_j}(x)| \, dx
\end{equation}
and note that 
$$
\nabla {\mathcal K}_{m^j} (x) = \int i \xi \phi(\xi) m(2^j \xi) e^{ i x \cdot \xi} \, d \xi 
\ =: \ \int \psi(\xi) m(2^j \xi) e^{ i x\cdot \xi} \, d\xi
$$
for some $\psi \in C^{\infty}_0({\mathbb R}^n)$ supported away from 0. Therefore 
$\nabla {\mathcal K}_{m^j}$ satisfies the bounds in \eqref{hypothesis-euclidean-pos}.
We write
$$
 \int_{|x| \ge |y|^{1-\theta}} |\nabla{\mathcal K}_{m_j}(x)| \, dx \ = \ 2^j
  \int_{|x| \ge 2^j |y|^{1-\theta}} |\nabla{\mathcal K}_{m^j}(x)| \, dx
$$ 
$$
 = \ 2^j \int_{2^j |y|^{1-\theta} \le |x| \le 2^{j\theta}} |\nabla{\mathcal K}_{m^j}(x)| \, dx \ + \
 2^j \int_{2^{j\theta} \le |x|} |\nabla{\mathcal K}_{m^j}(x)| \, dx \ =: \ I_j + II_j.
$$
We note that the integration in $I_j$ is nonempty since $|y|^{1-\theta} \le 2^{-j(1-\theta)}$
for $j\in J_2$. By Cauchy-Schwarz and \eqref{hypothesis-euclidean-pos} we have
$$
I_j \ \le \ 2^j 2^{j\theta n/2} \|{\mathcal K_{m^j}}\|_{L^2} \ = \ 
2^j 2^{j\theta n/2} \|m^j\|_{L^2} \ \lesssim \
2^j 2^{j\theta n/2} \|m^j\|_{L^{\infty}} \ \lesssim \ 2^j.
$$ 
In precisely the same way we argued in \eqref{same-way} we also have $|II_j| \lesssim 2^j$.
Hence $\sum_{j\in J_2} |I_j + II_j| \lesssim |y|^{-1}$ and this shows that we can sum the
integrals in \eqref{difference-y} and get a uniform bound, establishing the claim that \eqref{K-assumptions}
holds for ${\mathcal K}^l$.

\subsection{An interlude} At this point we would like to highlight a useful bound which is trivial in the euclidean setting
but will not be so trivial in the Lie group setting. The following
bound is an immediate consequence of the Cauchy-Schwarz inequality:
\vskip 15pt
\begin{mdframed}[style=JimFrame]
For any compactly support $F$ with ${\rm supp}(F) \subseteq K$ ($K$ compact),
\begin{equation}\label{key-euclidean}
\|{\mathcal K}_F \|_{L^1({\mathbb R}^n)} \, dx \ \le \ C_{s, K} \|F\|_{L^2_s({\mathbb R}^n)}
\end{equation}
holds for any $s > n/2$.
\end{mdframed}

We can use \eqref{key-euclidean} to conclude that if the decay parameter $\beta > n$, then
the main part of the convolution kernel ${\mathcal K}^l$ is integrable for any
$m\in {\mathcal M}_{\beta}$. To see this, note that $m \in M_{\theta, \beta, s}$ 
for some $\theta \in {\mathbb R}$ and $s > n/2$, and by \eqref{key-euclidean},
$$
\|{\mathcal K}^l \|_{L^1} \le \sum_{j: j\theta > 0} \| {\mathcal K}_{m_j}\|_{L^1}   = 
 \sum_{j: j\theta > 0} \| {\mathcal K}_{m^j}\|_{L^1}  \lesssim  \sum_{j: j\theta>0} 
\|m^j\|_{L^2_{s'}} \lesssim \sum_{j: j\theta > 0} 2^{ - j \theta (\beta - 2s')/2}
$$
for any $s' > n/2$. Since $\beta > n$ and $s>n/2$, we can find
an $s' \le s$ such that $n/2 < s' < \beta/2$.  Hence the above sum is convergent
and this shows that ${\mathcal K}^{l} \in L^1({\mathbb R}^n)$. 

By embedding a general $m \in {\mathcal M}_{\beta}$ with $0\le \beta < n$ into the analytic family of multipliers
$m_z (\xi) = |\xi|^{\theta/2(\beta - (n+\delta) z)} m(\xi)$ (see \eqref{analytic-family}) and using analytic interpolation,
we have the following observation.

\begin{lemma}\label{open-range} If $m\in {\mathcal M}_{\beta}$ and $0\le \beta < n$, then 
$m$ is an $L^p({\mathbb R}^n)$ multiplier if $|1/p - 1/2| < \beta/2n$.
\end{lemma}
Lemma \ref{open-range} is an extension of a result in \cite{S1} from the case $0<\theta < 1$ to the
case of general $\theta \not= 1$. 

\subsection{The results} As discussed above, using \eqref{M-beta} and the analytic interpolation argument in \cite{FS1},
we can show that any 
$m\in {\mathcal M}_{\beta}$ with $0\le \beta < n$ is an $L^p$ multiplier at the endpoint
$|1/p - 1/2| = \beta/2n$ IF we can show that every endpoint multiplier
$m \in {\mathcal M}_{n}$ is bounded from $H^1({\mathbb R}^n)$ to $L^1({\mathbb R}^n)$. 
We have the following theorem.

\begin{theorem}\label{main-euclidean} For every $m\in {\mathcal M}_{n}$,
the corresponding multiplier operator $T_m$ is weak-type $(1,1)$ and maps $H^1({\mathbb R}^n)$
to $L^1({\mathbb R}^n)$.
\end{theorem}
We do not claim that Theorem \ref{main-euclidean} is really new. For the examples
in \eqref{osc-multipliers}, Theorem \ref{main-euclidean} was established in the series
of papers \cite{F}, \cite{FS1}, \cite{CKW1}, \cite{CKW2} and \cite{M} for various
cases of $\theta \in {\mathbb R}\setminus\{1\}$. What is new is the proof which gives
a unified approach and extends to the Lie group setting. We have the immediate consequence 
improving Lemma \ref{open-range}.

\begin{corollary}\label{Lp-sharp}
If $m \in {\mathcal M}_{\beta}$ and $0\le \beta < n$, then $m$ is an $L^p({\mathbb R}^n)$ multiplier
for $|1/p - 1/2| \le 2\beta/n$.
\end{corollary}

\section{The stratified Lie group setting}\label{lie-group}

Let ${\mathfrak g}$ be an $n$-dimensional, graded nilpotent Lie algebra so that 
$$
{\mathfrak g} \ = \ \bigoplus\limits_{i=1}^s \, {\mathfrak g}_i  
$$
as a vector space and $[{\mathfrak g}_i, {\mathfrak g}_j] \subset {\mathfrak g}_{i+j}$ for
all $i,j$. Suppose that ${\mathfrak g}_1$ generates ${\mathfrak g}$ as a Lie algebra. We call 
the associated, connected, simply connected Lie group $G$ a stratified Lie group. Associated to
such a group is its so-called homogeneous dimension 
$$
Q \ = \ \sum_j j \, {\rm dimension}({\mathfrak g}_j)
$$
which is clearly always larger then the topological dimension $n$ but they agree when $G = {\mathbb R}^n$.

We fix
a basis $\{X_j\}$ for ${\mathfrak g}_1$ where each $X_j$ can be identified with a unique left-invariant
vector field on $G$ which we also denote by $X_j$.
Consider the sublaplacian ${\mathcal L} = - \sum_k X_k^2$ on $G$. For
any Borel measurable function $m$ on ${\mathbb R}_{+} = [0,\infty)$, 
we can define the spectral multiplier operator
$$
m(\sqrt{{\mathcal L}}) \ = \ \int_0^{\infty} m(\lambda) \, d E_{\lambda}
$$
where $\{E_{\lambda}\}_{\lambda\ge 0}$ is the spectral resolution of $\sqrt{{\mathcal L}}$. This is
a bounded operator on $L^2(G)$ precisely when $m \in L^{\infty}({\mathbb R}_{+})$.
The classical laplacian $\Delta$ is the corresponding differential operator when
$G = {\mathbb R}^n$ and spectral multipliers on ${\mathbb R}^n$ are simply radial
multipliers which the multipliers in \eqref{osc-multipliers} provide specific examples.

\subsection{The multiplier classes} 
We now state the conditions corresponding to \eqref{hypothesis-euclidean-neg} and
\eqref{hypothesis-euclidean-pos} for spectral multipliers $m$ defined
on ${\mathbb R}_{+}$. Fix a smooth bump function $\phi$ on ${\mathbb R}$ supported
away from the origin and let $m^j(\lambda) := m(2^j \lambda) \phi(\lambda)$. 
Again the conditions will depend on an oscillation parameter $\theta \in {\mathbb R}$,
a decay parameter $\beta \ge 0$ and a smoothness parameter $s>0$.
When $j \theta \le 0$, 
we impose the standard uniform $L^2$ Sobolev norm control on the $m_j$; 
\begin{equation}\label{hypothesis-lie-neg}
\sup_{j:  j \theta \le 0} \| m^j \|_{L^2_s({\mathbb R}_{+})} \ < \ \infty.
\end{equation}
For $j \theta > 0$, we consider the condition
\begin{equation}\label{hypothesis-lie-pos}
\sup_{j: j \theta >  0} \,  2^{j \theta \beta/2} \|m^j\|_{L^{\infty}({\mathbb R}_{+})},  
\  \ 2^{- j \theta (2s - \beta) /2} \|m^j\|_{L^2_s({\mathbb R}_{+})} \ < \ \infty. 
\end{equation}
Again when $\theta = 0$, these conditions reduce to the condition
$\sup_j \|m^j\|_{L^2_s} < \infty$ and if this holds for some $s > Q/2$, the fundamental work
of Christ \cite{C} and Mauceri-Meda \cite{MM} establishes that the multiplier operator is
of weak-type $(1,1)$ and bounded on  $H^1(G)$.

The examples 
$m_{\theta, \beta}(\lambda) =  e^{i \lambda^{\theta}} \lambda^{-\theta \beta/2} \chi_{\pm}(\lambda)$
from \eqref{osc-multipliers} satisfy conditions \eqref{hypothesis-lie-neg} and
\eqref{hypothesis-lie-pos}. As before we let $M_{\theta, \beta, s}$ consist of those spectral
functions $m$ satisfying \eqref{hypothesis-lie-neg} with exponent $s$ and satisfying
\eqref{hypothesis-lie-pos} for all exponents $s'\le s$. We redefine
\begin{equation}\label{M-beta-lie}
{\mathcal M}_{\beta} \ (= {\mathcal M}_{\beta, Q}) \ := \ \bigcup_{\theta \in {\mathbb R}\setminus\{1\}, s>Q/2} M_{\theta, \beta, s} \ = \Bigl\{
\lambda^{(Q-\beta)\theta/2} m(\lambda) :  m \in {\mathcal M}_Q \Bigr\}
\end{equation}
and stress the dependence of these classes on the homogeneous dimension $Q$ which we will return to later.
Again this puts us in the position to employ analytic interpolation arguments 
to deduce that $m \in {\mathcal M}_{\beta}$ is an $L^p(G)$ multiplier in the range $|1/p - 1/2| \le \beta/2Q$
from $H^1(G)$ bounds for multiplier operators associated to $m\in {\mathcal M}_Q$. 
Furthermore, from the invariance of ${\mathcal M}_Q$ under multiplication by $\lambda^{iy}$
for any real $y$ (with resulting polynomial in y bounds in \eqref{hypothesis-lie-neg} and \eqref{hypothesis-lie-pos}),
it suffices to show $m(\sqrt{{\mathcal L}}) : H^1(G) \to L^1(G)$ for $m\in {\mathcal M}_Q$.

\subsection{The main result} Our main result is the following theorem.

\begin{theorem}\label{main-lie} For any $m\in {\mathcal M}_Q$, the operator 
$m({\sqrt{{\mathcal L}}}) : H^1(G) \to L^1(G)$ and is weak-type $(1,1)$. 
\end{theorem} 
As an immediate consequence, using analytic interpolation (see above), we have the following
endpoint result of Mauceri and Meda in \cite{MM}. See also the work of Alexopolous \cite{Alex}
on general Lie groups of polynomial volume growith.

\begin{corollary}\label{corollary-lie} Every $m\in {\mathcal M}_{\beta}$ with $0\le \beta < Q$
is an $L^p(G)$ multiplier in the range $|1/p - 1/2| \le \beta/2Q$. 
\end{corollary}

\subsection{The interlude - revisited} We now return to the estimate \eqref{key-euclidean}
and examine it in the Lie group context. Again we use the notation ${\mathcal K}_F$ to denote
the convolution kernel of the operator $m(\sqrt{{\mathcal L}})$.
\vskip 15pt
\begin{mdframed}[style=JimFrame]
Let $G$ be any stratified Lie group and suppose the following
holds for some dimensional parameter $d$: for
any spectral multiplier $F(\lambda)$, supported in a compact $K\subset {\mathbb R}_{+}$, 
\begin{equation}\label{key-lie}
\|{\mathcal K}_F \|_{L^1(G)} \, dx \ \le \ C_{s, K} \|F\|_{L^2_s({\mathbb R}_{+})}
\end{equation}
holds for any $s > d/2$.
\end{mdframed}

In \cite{C} and \cite{MM}, the estimate \eqref{key-lie} was proved for $d=Q$, the homogeneous dimension,
on a general stratified Lie group $G$. In fact the estimate \eqref{key-lie} is the key estimate in
their work. It is known that if \eqref{key-lie} holds for some parameter $d$, then standard techniques
allow us to deduce that if a spectral multiplier $m$ satsifies $\sup_j \|m^j\|_{L^2_s} < \infty$
for some $s> d/2$, then $m(\sqrt{{\mathcal L}})$ is bounded on all $L^p(G), 1< p < \infty$
and corresponding endpoint results on $L^1$ hold. See for example, \cite{alessio-joint}. Hence
to determine the minimal amount of smoothness required
for H\"ormander-type spectral multipliers, matters can be reduced to establishing \eqref{key-lie}.

The fact that one only needs to control a little more than half the topological dimension $n$
number of derivatives, $s > n/2$, for certain Lie groups was first observed by M\"uller and Stein
\cite{MS} for the Heisenberg group. The ideas in \cite{hebisch} can be used to establish
\eqref{key-lie} for $d=n$ on any Lie group of Heisenberg-type (alternatively, one of the main estimates in
\cite{MSeeger} imply this immediately). Furthermore \eqref{key-lie} for $d=n$ was established
by Martini and M\"uller \cite{plms} for step 2 stratified Lie groups with $n\le 7$ or whose
centre has dimension at most 2. In another paper \cite{gafa}, Martini and M\"uller show
that \eqref{key-lie} holds for some $d<Q$ on any step 2 stratified Lie group.

The estimate \eqref{key-lie} also has implications for our more general multipliers satsifying
\eqref{hypothesis-lie-neg} and \eqref{hypothesis-lie-pos}. Instead of 
${\mathcal M}_{\beta} = {\mathcal M}_{\beta, Q}$ defined in \eqref{M-beta-lie}, let us consider
$$
{\mathcal M}_{\beta, d} \ := \ \bigcup_{\theta \in {\mathbb R}\setminus\{1\}, s>d/2} M_{\theta, \beta, s} 
$$
depending now on a dimensional paramter $d$ which could be smaller than $Q$.
Suppose now that \eqref{key-lie} holds for some $d\le Q$ on $G$. 
We can use \eqref{key-lie} to conclude that if $\beta > d$, then any 
$m\in {\mathcal M}_{\beta, d}$ can be written as $m = m_{small} + m_{large}$ (see \eqref{m-split})
where $m_{small}$ is a H\"ormander multiplier with $s> d/2$ (and hence bounded on all
$L^p(G), 1<p<\infty$, weak-type $(1,1)$, etc...) and $m_{large}$ is an $L^1(G)$ multiplier,
the convolution kernel ${\mathcal K}^l$ associated to $m_{large}$ being integrable. This follows
exactly as in the the euclidean setting.

By embedding a general $m \in {\mathcal M}_{\beta, d}$ with $0\le \beta < d$ into the analytic family of 
spectral multipliers
$m_z (\lambda) = \lambda^{\theta/2(\beta - (d+\delta) z)} m(\lambda)$ and using analytic interpolation,
we have the following observation.

\begin{lemma}\label{open-range-lie} Suppose that \eqref{key-lie} holds on $G$ for some 
$d\le Q$. 
If $m\in {\mathcal M}_{\beta,d}$ and $0\le \beta < d$, then 
$m$ is an $L^p(G)$ multiplier for $|1/p - 1/2| < \beta/2d$.
\end{lemma}
In particular on any step 2 stratified Lie group, the result of Martini and M\"uller in \cite{gafa},
establishing that \eqref{key-lie} holds for some $d<Q$, shows that the convolution kernel ${\mathcal K}^l$
corresponding to the interesting frequency range of any $m \in {\mathcal M}_Q = {\mathcal M}_{Q,Q}$
is integrable! In this case our main result Theorem \ref{main-lie} is a consequence of the work of
Christ \cite{C} and Mauceri-Meda \cite{MM}.

Hence we should view Theorem \ref{main-lie} and Corollary \ref{corollary-lie} as place-holders
for possible endpont results. It may be the case that \eqref{key-lie} holds for some $d<Q$ on
any stratified Lie group outwith the euclidean $G = {\mathbb R}^n$ case. If so,  our results do
not say anything new outside the euclidean setting.

In a forthcoming paper, we will establish the sharp result on any Lie group of Heisenberg-type,
establishing Theorem \ref{main-lie} and Corollary \ref{corollary-lie} with $Q$ replaced by $n$. Our
analysis heavily relies on M\"uller and Seeger's work \cite{MSeeger} on the wave equation in Lie groups
of Heisenberg-type.

Finally we note that Theorem \ref{main-lie} implies Theorem \ref{main-euclidean} in the case
of radial multipliers but the proof of Theorem \ref{main-lie} below easily gives
a proof of Theorem \ref{main-euclidean}. We will therefore give the proof
of Theorem \ref{main-lie} only.

\section{Preliminaries}\label{preliminaries}

For background information about Calder\'on-Zygmund theory and spectral multipliers
on stratified groups, we refer the reader to the book of Folland and Stein \cite{FS-book}. 
If $h$ is a Borel measurable function on ${\mathbb R}_{+}$, recall that 
${\mathcal K}_h$ denotes the convolution kernel of the
operator $h(\sqrt{\mathcal L})$ so that
$$
h({\sqrt{\mathcal L}}) f (x) \ = \ f * {\mathcal K}_h (x) \ = \ \int_G f(x\cdot y^{-1}) \,
{\mathcal K}_h (y) \, dy
$$
where $dy$ denotes Haar measure on $G$. Since we are identifying the Lie group $G$
with its Lie algebra ${\mathfrak g}$ via the exponential map, the Haar measure is identified with
Lebesgue measure on the Lie algebra ${\mathfrak g} \simeq {\mathbb R}^n$.

\subsection{Some basics} The stratified group $G$ comes equipped with
a group of dilations $\delta_r : G \to G$ which are automorphisms and we fix 
a homogeneous norm; that is, a function $|\cdot|: G \to {\mathbb R}_{+}$, smooth 
away from $0$, with $|x| = 0$ if and only if $x=0$ where $0$ denotes the group identity, 
and $|\delta_r x | = r |x|$ for all $r\in {\mathbb R}_{+}$ and $x\in G$.
Also if $s>0$, then
$$
h(s \sqrt{\mathcal L}) f (x) = f * ({{\mathcal K}_h})_{s} (x) \ \ {\rm where} \ \
 ({{\mathcal K}_h})_{s} (x) \ := \  s^{-Q} {\mathcal K}_h (\delta_{s^{-1}} x);
$$
see \cite{FS-book}. Another standard fact from \cite{FS-book} is the following mean value
theorem for Schwartz functions ${\mathcal S}$ on $G$: if $h \in {\mathcal S}(G)$,
then for any $N\ge 1$, 
\begin{equation}\label{mean-value}
|h(x \cdot y) - h(x)| \ \le \ C_N \, \frac{|y|}{(1 + |x|)^N}
\end{equation}
holds for any $y\in G$ such that $|y| \ll |x|$. We will find this useful at times. We will
also find useful the following Plancherel-type identity which
can be found in \cite{C}: for $h \in L^{\infty}({\mathbb R}_{+})$, there is a constant
$c$ such that
\begin{equation}\label{plancherel}
\|{\mathcal K}_h\|_{L^2(G)}^2 \ = \ c \int_0^{\infty} |h(t)|^2 \ t^{Q -1} \, dt
\end{equation}
holds. 

\subsection{A weighted $L^2$ bound} We will use the following weighted $L^2$ estimate which is valid on a general
stratified Lie group $G$: if $F$ is a compactly supported spectral
multiplier, then
\begin{equation}\label{sikora}
\int_G |{\mathcal K}_{F}(x)|^2 (1+|x|^s)^2 \, dx \ \lesssim \ \|F\|_{L^2_s}
\end{equation}
holds for any $s > 0$. See \cite{sikora}. Note that by Cauchy-Schwarz, the bound 
\eqref{sikora} immediately shows that the key estimate \eqref{key-lie} holds for all $s>Q/2$ on any stratified Lie group.

For the Hardy space estimate we will use \eqref{sikora} but we will also use this estimate with derivatives:
\begin{equation}\label{sikora-derivatives}
\int_G |X_j K_{F({\mathcal L})}(x)|^2 (1+|x|^s)^2 \, dx \ \lesssim \ \|F\|_{L^2_s}
\end{equation}
holds for any $s>0$, $1\le j \le k$ and any compactly supported $F$. Here $k = {\rm dimension}({\mathfrak g}_1)$.

\subsection{Fefferman-Stein inequality} Our argument uses the Fefferman-Stein vector-valued Hardy-Littlewood maximal function
inequality in the context of stratified groups. If
$$
M f(x) \ = \ \sup_{r>0} \frac{1}{r^Q} \int_{|y|\le r} |f(x \cdot y^{-1})| \, dy 
$$
denotes the Hardy-Littlewood maximal function on $G$, then for
$1< p, q < \infty$, we have
\begin{equation}\label{FS-inequality}
\Bigl\| \Bigl( \sum_j (M f_j)^q \Bigr)^{1/q} \Bigr\|_{L^p(G)} \ \le \ C_{p,q,G} \,
\Bigl\| \Bigl( \sum_j |f_j|^q \Bigr)^{1/q} \Bigr\|_{L^p(G)};
\end{equation}
see for example \cite{S} or \cite{GLY}. 
We will use this inequality for $f_j$ a sequence of characteristic functions of balls
$B = B(x_B, r_B) := \{y \in G : |y\cdot x_B^{-1}| \le r_B\}$. We first note that if
$\chi_B$ denotes the characteristic function of a ball $B$, then
\begin{equation}\label{max-char}
M(\chi_B)(x) \ \sim \ \frac{1}{(1 + |\delta_{2^{-L(B)}}( x \cdot x_B^{-1})| )^Q}
\end{equation}
where $L(B)$ is chosen so that $2^{L(B)} = r_B$.
Hence $M(\chi_B)$ is a weak approximation
of the characteristic function $\chi_B$ itself. 

\subsection{Our basic decomposition} Let us recall the basic decomposition \eqref{m-split} in the context of spectral multipliers $m$; we
choose $\phi \in C^{\infty}_0({\mathbb R}_{+})$ supported in $\{ 1/2 \le \lambda \le 2\}$
so that $\sum_{j\in {\mathbb Z}} \phi(2^{-j} \lambda) = 1$ for all $\lambda>0$. 
Hence $m(\lambda) = \sum_{j\in {\mathbb Z}} m_j (\lambda)$ where 
$m_j(t) := m(t) \phi(2^{-j} \lambda) = m^j (2^{-j} \lambda)$ and so  
\begin{equation}\label{mj}
{\mathcal K}_{m_j} (x)  \ = \ {\mathcal K}_m  * ( 2^{jQ} {\mathcal K}_{\phi} (\delta_{2^{j}} \cdot) ) (x)
\ = \ {\mathcal K}_m * ({\mathcal K}_{\phi})_{2^{-j}} (x).
\end{equation}
Therefore
$$
m(\sqrt{\mathcal L}) f (x)\  = \ \sum_{j\in {\mathbb Z}} m_j(\sqrt{\mathcal L}) f(x) \ = \ 
\sum_{j\in {\mathbb Z}} f * {\mathcal K}_{m_j} (x).
$$

For $m \in {\mathcal M}_Q$, we split the multiplier
$$
m(\lambda) \ = \ \sum_{j\in {\mathbb Z}} m_j (\lambda)  \ = \sum_{j\in {\mathbb Z}} \ 
m(\lambda) \phi(2^{-j}\lambda) \ = \ m_{small}(\lambda) + m_{large}(\lambda)
$$
into two parts 
where $m_{small}(\lambda) = \sum_{j: j \theta \le 0} m_j (\lambda)$ and 
$m_{large}(\lambda) = \sum_{j: j\theta > 0}  m_j (\lambda)$.
Since $m$ satisfies \eqref{hypothesis-lie-neg} for some $s>Q/2$, the results of Christ
\cite{C} and Mauceri-Meda \cite{MM} show the multiplier $m_{small}$ is weak-type $(1,1)$
and bounded on $H^1(G)$ (alternatively, the argument below in the case $\theta = 0$
can be used to treat $m_{small}$). Hence it suffices to  treat the operator
 $T := \sum_{j : j \theta > 0} m_j (\sqrt{\mathcal L})$
 and in particular it will
be good to keep in mind that $j \theta > 0$ is always satisfied.

\section{The proof of Theorem \ref{main-lie} -- the weak-type $(1,1)$ bound}

We have reduced matters to bounding $T = \sum_{j : j \theta > 0} m_j (\sqrt{\mathcal L})$
and our aim here is to show that
\begin{equation}\label{1-1}
\Bigl|\bigl\{ x \in G:  |T f(x)| \ge \alpha \bigr\}\Bigr| \ \le \  \frac{C}{\alpha}  \,
\|f\|_{L^1(G)}
\end{equation}
holds uniformly for all $\alpha > 0$ and $f \in L^1(G)$. We will denote by $|\cdot |$ the
Haar measure on $G$ as well as the homogeneous norm (as well as the usual absolute
value on ${\mathbb R}$ or ${\mathbb C}$). There should be no confusion.

We employ the classical Calder\'on-Zygmund
decomposition of $f$ at height $\alpha$ on $G$ (see \cite{FS-book}
or \cite{S}): 
there exists a sequence of essentially disjoint balls
$\{ B = B(x_B, 2^{L(B)}) \}$ such that $|\cup B| \lesssim \|f\|_{L^1}/\alpha$. Furthermore
we can decompose
$f = g + b$ where $|g(x)| \lesssim \alpha$ a.e $x\in G$
and $b = \sum_B b_B$ where ${\rm supp}(b_B) \subseteq B^{*}$, 
\begin{equation}\label{b-properties}
\int_G b_B = 0, \ \ \|b_B\|_{L^1} \lesssim \alpha |B| \ \ \ {\rm and} \ \
\sum_B \|b_B \|_{L^1} \lesssim \|f\|_{L^1}.
\end{equation}
Here and from now on, $B^{*}$ will
denote a generic dilate of $B$ which is understood to be the appropriate dilate depending on
the context and we may also take it to be a sufficiently large dilate when there is a need to do so.

The contribution of the bounded function $g$ to the distribution function 
$|\{ x : |Tf(x)| \ge \alpha \}|$ follows
in the usual way, only the $L^2$ boundedness of $T$  is used here 
(that is, only the fact that $m$ is bounded is used). 
To establish \eqref{1-1}, it suffices to
consider the contribution from $T$ on the function $b = \sum_B b_B$ where
$f$ is large and so we write
$$
T b(x) \ = \  \sum_{(j,B)\in N} {\mathcal K}_{m_j} * b_B (x) \ + \ \sum_{(j,B)\in P} 
{\mathcal K}_{m_j} * b_B(x) \ =: \ {\mathcal A}(x) + {\mathcal B}(x)
$$
where 
$$
N = \{(j,B) : j \theta > 0, \  j (1-\theta) +  L(B) \le 0  \}
$$
and $N$ is the complementary set of pairs $(j,B)$ with $j \theta > 0$. 

For the sum over the pairs $(j,B)\in N$, we use $L^2$ estimates,
the disjoint frequency supports of the $\{\phi(2^{-j}\lambda)\}$ and 
the smallness of $m$
on the support of $\phi(2^{-j}\lambda)$, $m \approx 2^{-\theta j Q/2}$. Writing
$\Phi_j(x) := ({\mathcal K}_{\phi})_{2^{j}}(x)$, we have
$$
|\{ x: |\sum_{(j,B)\in N} m_j(\sqrt{\mathcal L}) (b_B)(x) | \ge \alpha \}| \le \alpha^{-2} \| 
\sum_{(j,B)\in N} m_j(\sqrt{\mathcal L})( b_B)\|_2^2 
$$
$$
\lesssim \alpha^{-2} \sum_{j: j \theta > 0} \| {\mathcal K}_{m} * 
\bigl(\sum_{B \in  N_j } \Phi_{j} * b_B \|_2^2 \ \lesssim \
\alpha^{-2} \sum_{j: j \theta > 0} 2^{-\theta j Q} \| \sum_{B\in N_j} \Phi_{j} * b_B \|_2^2
$$
where $N_j = \{ B : (j,B)\in N\}$.
We write the last term on the right above as $E + F$
where
$$
E := \alpha^{-2} \sum_{j: j \theta > 0} 2^{- j Q\theta} \| 
\sum_{B\in N_j} \Phi_{j} * b_B \cdot \chi_{B^{*}} \|_2^2
\lesssim \alpha^{-2} \sum_{(j,B)\in N} 2^{- j Q\theta} \| \Phi_{j} * b_B  \|_2^2
$$
for some appropriately large dilate $B^{*}$ of $B$
and $F$ is defined similarly with $B^{*}$ replaced by 
$G\setminus B^{*}$.
Since $\| \Phi_{j} * b_B\|_{L^2}^2 \le \| \Phi_{j}  \|_{L^2}^2 \|b_B\|_{L^1}^2
\lesssim \|\Phi_j \|_{L^2(G)}^2\alpha^2 |B|^2$
and
$$
\|\Phi_j \|_{L^2(G)}^2 \ = \ c \int_0^{\infty} |\phi(2^{-j} \lambda)|^2 t^{Q -1} \, dt
\ = \ c_{\phi} 2^{j Q}
$$
by the Plancherel formula \eqref{plancherel}, we have
$ \| \Phi_{j} * b_B\|_{L^2}^2 \lesssim 2^{j Q} \alpha^2 |B|^2$.
Hence 
$$
E \lesssim \sum_{(j,B)\in N} 2^{j Q(1-\theta)} |B|^2 \lesssim \sum_B |B|
\sum_{j(1-\theta)+ L(B) \le 0 } 
2^{Q(j(1-\theta) + L(B))}
\lesssim \sum_B |B| \lesssim \alpha^{-1} \|f\|_1.
$$
Note that it is important that $\theta \not=1$ in the above argument. This
leaves us with $F$. 

Using the cancellation of $b_B$, we have
$$
\Phi_{j}* b_B(x) \ = \
\int_G \bigl[ \Phi_{j}(x\cdot y^{-1}) - \Phi_{j}(x\cdot x_B^{-1}) \bigr] b_B(y) \, dy .
$$
Noting that $\Phi_j (x) = 2^{j Q} {\mathcal K}_{\phi} (\delta_{2^{j}} x)$, we have 
for $y \in {\rm supp}(b_B)$ and $x \notin B^{*}$,
$$
|\Phi_j (x \cdot y^{-1}) - \Phi_j (x \cdot x_B^{-1}) | \ \lesssim \ 2^{j Q} \,
\frac{2^{(1-N)(j + L(B))}}{(1 + |\delta_{2^{-L(B)}}( x \cdot x_B^{-1}) |)^N}
$$
by the mean value theorem on stratified groups \eqref{mean-value}. Therefore we see that for $x\notin B^{*}$,
$$
|\Phi_j * b_B (x) | \lesssim \alpha 2^{(Q + 1-N)(j + L(B))} \, M(\chi_B)(x)^{N/Q} = 
\alpha 2^{\epsilon (j + L(B))} M(\chi_B)(x)^q
$$
where $\epsilon = Q+1-N$ and $q = N/Q$. By choosing
$N = Q + 1/2$, we can make $\epsilon > 0$ and $q>1$. This allows us to apply
the Fefferman-Stein inequality \eqref{FS-inequality} which yields
$$
F \ \lesssim \ \alpha^{-2} \sum_{j: j \theta > 0} 2^{- j Q\theta}  
\bigl\| \sum_{B\in N_j} \Phi_{j} * b_B (\chi_{G \setminus B^{*}}) \bigr\|_2^2
$$
$$
\lesssim \sum_{j: j\theta > 0} 2^{- j Q\theta}  \bigl\| 
\sum_{B\in N_j} [M(2^{\epsilon(j+L(B))/q}\chi_B)]^q \|_2^2
\lesssim \sum_{(j,B)\in N} 2^{- j Q\theta} 2^{2\epsilon (j+L(B))} |B| 
$$
$$
= \sum_{(j,B)\in N} 2^{j(1-\theta) + L(B)} 
2^{-j\theta (Q-1)} |B| 
\le \sum_{(j,B)\in N} 2^{j(1-\theta) + L(B)}  |B| \lesssim  \sum_B |B| \ \lesssim \ \alpha^{-1} \|f\|_1
$$
since $j\theta > 0$ and $Q\ge 1$. 

This completes the estimate for $F$ and the contribution from the pairs $(j,B) \in N$.
Hence
$|\{x: |{\mathcal A}(x)| \ge \alpha \}| \lesssim \alpha^{-1} \|f\|_1$.
Again it was important
that $\theta\not=1$ in this argument. We now turn to
the contribution from the pairs $(j,B)\in P$ where we will use only $L^1$ estimates and the $L^2$ 
Sobolev condition in \eqref{hypothesis-lie-pos}.

Since $|\cup B^{*}| \lesssim \alpha^{-1} \|f\|_1$, we see that the desired estimate
$|\{x : |{\mathcal B}(x)| \ge \alpha\}| \lesssim \alpha^{-1} \|f\|_1$ reduces matters to 
estimating $|\{ x \notin \cup B^{*} : |{\mathcal B}(x) | \ge \alpha \}|$ which we see is at most
$$
 \alpha^{-1} 
\int_{x\notin \cup B^{*}} |\sum_{(j,B)\in P} m_j(\sqrt{\mathcal L})(b_B)(x)| dx \ \le \ \alpha^{-1}
\sum_{(j,B)\in P} \int_{x\notin B^{*}}  |m_j(\sqrt{\mathcal L}) (b_B)(x)| dx 
$$
$$
\le \ 
\alpha^{-1} \sum_{(j,B)\in P} 
\int |b_B(y)| \Bigl[\int_{|x\cdot x_B^{-1}|\gg 2^{L(B)}} |{\mathcal K}_{m_j}(x\cdot y^{-1})| dx \Bigr] dy .
$$
The desired estimate will follow
if we can show that 
\begin{equation}\label{desire}
\sup_B \, 
\sum_{j: j(1-\theta) + L(B) \ge 0} \int_{|x|\gtrsim 2^{L(B)}} |{\mathcal K}_{m_j}(x)|\,  dx \ < \ \infty .
\end{equation}

In fact,
$$
\int_{|x|\ge 2^{L(B)}} |{\mathcal K}_{m_j}(x)| dx \ = \ 
 \int_{|x|\ge 2^{j + L(B)}} |{\mathcal K}_{{m}^j}(x) | dx
$$
$$
2^{(Q/2 - s)(j + L(B) )} \int_G |K_{{m}^j}(x) |^2 (1 + |x|^s)^2 dx \ \lesssim \
  2^{-(s-Q/2)(j(1-\theta) + L(B))} 
$$
and this sums in $j$ with $j(1-\theta) + L(B) \ge 0$ if $s>Q/2$, uniformly in $B$. 
Here we used \eqref{sikora} and the $L^2$ Sobolev condition in \eqref{hypothesis-lie-pos} in the penultimate inequality. This establishes \eqref{desire} and completes the proof of the weak-type $(1,1)$ bound in Theorem \ref{main-lie}.

\section{The proof of Theorem \ref{main-lie} -- the Hardy space bound}

Elements in the Hardy space $H^1(G)$ have an atomic deomposition (see \cite{FS-book}) and so
it suffices to fix an atom 
 $a_B$ supported in a ball $B$ and prove
\begin{equation}\label{L1-atom}
\int_G |m(\sqrt{\mathcal L}) a_B (x) | \, dx \ \lesssim \ 1
\end{equation}
for our spectral multiplier $m \in {\mathcal M}_Q$. 

Without loss of generality we may assume that the ball $B$ is centred at the origin. The $L^2$ boundedness of $m(\sqrt{\mathcal L})$
implies that $\int_{|x|\le C 2^L} |m(\sqrt{\mathcal L}) a_B(x)| dx \lesssim 1$ via the Cauchy-Schwarz inequality and so it suffices
to show that 
\begin{equation}\label{desire-2}
\int_{|x| \gg 2^L} |m(\sqrt{\mathcal L}) a_B (x)| \, dx \ \lesssim \ 1
\end{equation}
holds where $2^L$ is the radius of the ball $B$.

From our basic decomposition $m = m_{small} + m_{large}$, it suffices as before to treat
the operator 
$T := \sum_{j : j \theta > 0} m_j (\sqrt{\mathcal L})$ and show that \eqref{desire-2} holds
with $m(\sqrt{\mathcal L})$ replaced by $T$.

We bound the integral in \eqref{desire-2} by
$$
\sum_{j\in N} \int_{|x| \gg 2^L} |{\mathcal K}_{m_j} * a_B (x)| dx \ + \
\sum_{j \in P} \int_{|x| \gg 2^L} |{\mathcal K}_{m_j} * a_B (x)| dx =: I + II
$$
where $N = \{ j: j \theta > 0,  j(1-\theta) + L \le 0 \}$ and $P$ denotes the complementary range. 

For $j\in P$, we note that when $|x|\gg 2^L$,
$$
{\mathcal K}_{m_j} * a_B (x) \ = \int_{|y|\le 2^L} {\mathcal K}_{m_j}(x y^{-1}) a_B (y) \, dy \ = \
\int_G {\mathcal K}_{m_j}(x y^{-1}) \chi_{E_L}(x y^{-1}) a_B(y) \, dy
$$
where $E_L = \{ x \in G : |x| \ge 2^L\}$. Hence if we denote by $K = {\mathcal K}_{m_j}(x) \chi_{E_L}(x)$,
$$
\int_{|x|\gg 2^L} |{\mathcal K}_{m_j} * a_B(x)| dx = \int_G |K * a_B(x) | dx \le \int_G |K(x)| dx =
\int_{|x|\ge 2^L} |{\mathcal K}_{m_j}(x)| dx
$$
$$
= \ \int_{|x| \ge 2^{j + L}} |{\mathcal K}_{{m}^j}(x)| \, dx \ = \ 
\int_{|x| \ge 2^{j + L }} |{\mathcal K}_{{m}^j}(x)| \, \frac{1+|x|^s}{1+|x|^s} \, dx
$$
$$
\le \ 2^{-(s - Q/2)(j + L)} \sqrt{\int_G |{\mathcal K}_{m^j}(x)|^2 (1+ |x|^s)^2 dx}
\ \lesssim \ 2^{-(s-Q/2)(j(1-\theta) + L)}
$$
where in the last inequality we used \eqref{sikora} with some $s> Q/2$ and the $L^2$ 
Sobolev condition of our multiplier $m$ as stated in
\eqref{hypothesis-lie-pos}. Since $\theta \not= 1$, this shows that $II$ is uniformly bounded

For $I$, we split $N = N_1 \cup N_2$ further such that $N_1 = \{ j \in N: j + Q \le 0\}$ and
$N_2 = \{j \in N: j + Q > 0\}$. This splits $I = I_1 + I_2$ accordingly.

For the sum over $j\in N_1$, we will use the cancellation of
the atom $a_B$: for $j\in N_1$,
$$
 \int_{|x| \gg 2^L} |{\mathcal K}_{m_j} * a_B (x)| dx \le \int_G |a_B(y)| \Bigl[\int_{C2^{j +L} \le |x|}
 |{\mathcal K}_{{m}^j } (x (\delta_{2^{j}} y)^{-1}) - {\mathcal K}_{{m}^j} (x)| dx \Bigr] dy
 $$ 
 and so by applying the mean value theorem on stratified groups (see \eqref{mean-value}), 
 we see that the inner integral on the right hand side is at most
 $$
 2^{j + L} \ \int_{2^{j+L} \le |x|} \sup_{1\le r\le k} |X_r {\mathcal K}_{{m}^j}  (x)| \, dx
 $$
 and so
 $$
 \int_{|x| \gg 2^L} |{\mathcal K}_{m_j} * a_B (x)| dx \ \le \ 2^{j +L} \ \sum_{r=1}^k  \
 \int_{2^{j+L} \le |x|} |X_r {\mathcal K}_{{m}^j}  (x)| \, dx.
 $$
Let $X$ denote one of the $X_r$'s -- our immediate goal is to show that the bound
\begin{equation}\label{desire-3}
\int_{2^{j+L} \le |x|} |X {\mathcal K}_{{m}^j} (x)| \, dx \ \le \ C
 \end{equation}
 holds, uniformly for all $j$ and $L$. If this is the case, then we see that
 $$
I_1 \ =  \ \sum_{j\in N_1} \int_{|x| \gg 2^L} |{\mathcal K}_{m_j} * a_B (x)| \, dx \ \lesssim \ \sum_{j\in N_1} 2^{j + L} 
\ \lesssim \ 1,
$$
completing the analysis for $I_1$. 

To show \eqref{desire-3}, we will use \eqref{sikora-derivatives}
for two different values of $s$. We split the integral in \eqref{desire-3} into two parts:
$$
\int_{2^{j+L} \le |x| \le 2^{j + L + \Lambda}} |X {\mathcal K}_{{m}^j } (x)| \, dx \ + \ 
\int_{2^{j+L + \Lambda} \le |x|} |X {\mathcal K}_{{m}^j} (x)| \, dx  =: S_{\Lambda} + L_{\Lambda}
$$
for some large $\Lambda>0$ to be chosen appropriately.

For $S_{\Lambda}$ we use \eqref{sikora-derivatives} with some $s_{*} < Q/2$: by Cauchy-Schwarz,
$$
S_{\Lambda}^2  \ \le\  2^{2(Q/2 - s_{*})(j + L + \Lambda)} \int_G  
|X {\mathcal K}_{{m}^j } (x)|^2 (1+|x|^{s_{*}})^2 \, dx
$$
and so using \eqref{sikora-derivatives} and the $L^2$ Sobolev condition 
\eqref{hypothesis-lie-pos} of our multiplier $m$, 
$$
S_{\Lambda} \ \lesssim \ 2^{(Q/2 - s_{*})( j (1 - \theta) + L)} \, 2^{(Q/2 - s_{*}) \Lambda}.
$$
In a similar way, using \eqref{sikora-derivatives} with some $s> Q/2$, we have
$$
L_{\Lambda}^2  \ \le\  2^{-2(s - Q/2)(j + L + \Lambda)} \int_G  |X {\mathcal K}_{{m}^j} (x)|^2 
(1+|x|^{s})^2 \, dx
$$
and so using the $L^2$ Sobolev condition \eqref{hypothesis-lie-pos} of our muliplier $m$, we see that
$$
L_{\Lambda} \ \lesssim \ 2^{-(s - Q/2)( j (1 - \theta) + L)} \, 2^{-(s - Q/2) \Lambda}.
$$
Optimising the two estimates gives $\Lambda = - (j(1-\theta) + L)$ which is positive since $j\in N$. Hence with
this choice of $\Lambda$,
$S_{\Lambda} + L_{\Lambda} \lesssim 1$, establishing \eqref{desire-3} and completing the analysis for $I_1$.

Finally we turn to $I_2$ where $j\in N_2$ implies $j + L \ge 0$. Here it does not make sense to use the cancellation
of the atom $a_B$. Instead we use our knowledge of the $L^2$ size of $a_B$; $\|a_B\|_{L^2(G)} \le |B|^{-1/2} = 2^{-LQ/2}$.
We begin by splitting the integral into two parts as above:
$$
\int_{2^L \le |x| \le 2^{L+\Lambda}} |{\mathcal K}_{m_j} * a_B(x)| \, dx \ + \
\int_{|x|\ge 2^{L+\Lambda}} |{\mathcal K}_{m_j} * a_B(x)| \, dx \  : = \ S_{\Lambda} + L_{\Lambda}
$$
for some appropriate $\Lambda$. For $S_{\Lambda}$, we use the $L^{\infty}$ condition in
\eqref{hypothesis-lie-pos} and Cauchy-Schwarz to see that
$$
S_{\Lambda} \le 2^{(L+\Lambda)Q/2} \|{\mathcal K}_{m_j} * a_B\|_{L^2} \le 
2^{(L+\Lambda)Q/2} 2^{-j\theta Q}
\|\Phi_j * a_B \|_{L^2} \le  2^{\Lambda Q/2} 2^{-j\theta Q/2}.
$$
On the other hand, for $L_{\Lambda}$, we have 
$$
L_{\Lambda} \ \le \ \int_{|x|\ge 2^{L+\Lambda}} |{\mathcal K}_{m_j}(x)| \, dx  \ = \
\int_{|x|\ge 2^{j + L+\Lambda}} |{\mathcal K}_{{m}^j}(x)| \, dx
$$
$$
\le \ 2^{-(s-Q/2)(j + L + \Lambda)} \sqrt{ \int_G |{\mathcal K}_{{m}^j}(x)|^2 (1 + |x|^s)^2 \, dx} \ \lesssim \
2^{-(s-Q/2)(j(1-\theta) + L + \Lambda)}
$$
by \eqref{sikora-derivatives} with $s>Q/2$ and Cauchy-Schwarz. Optimising the two estimates gives 
$\Lambda$ with $2^{s \Lambda} = 2^{-(s-Q/2) (j + L) } 2^{j \theta s /2}$ which is positive since $j\in N$. Hence with
this choice of $\Lambda$,
$S_{\Lambda} + L_{\Lambda} \lesssim 2^{-Q/2(1-Q/2s)/ (j + L)}$ which is summable over
$j\in N_2$ since $j + L > 0$,  showing that
$$
I_2 \ = \ \sum_{j\in N_2} \int_{2^L \le |x| } |{\mathcal K}_{m_j} * a_B(x)| \, dx
$$
is uniformly bounded in $L$ and this completes the analysis for $I_2$, establishing \eqref{desire-2}
and hence \eqref{L1-atom}. 

This finishes the $H^1(G)$ bound of $m(\sqrt{\mathcal L})$ and hence the proof of Theorem \ref{main-lie}.

\end{document}